\documentclass[11pt,a4paper,oneside,reqno]{amsart}

\usepackage{amsmath,amsthm,amssymb,amsopn}
\usepackage{hyperref}
\usepackage[capitalise]{cleveref}
\usepackage{enumerate}
\usepackage{bbm}
\usepackage{stmaryrd}
\usepackage[a4paper,total={17cm,25cm}]{geometry}
\usepackage{thmtools}
\usepackage{mathtools}
\usepackage{graphicx}
\usepackage{caption}
\usepackage{subcaption}
\usepackage{float}
\usepackage[font=small,labelfont=bf]{caption}
\usepackage{setspace}

\usepackage{xcolor}

\declaretheorem[name=Theorem]{thm}

\newcommand{\prodG}{\mathbf{G}}

\newcommand{\br}[1]{\llbracket{#1}\rrbracket}

\renewcommand{\le}{\leqslant}
\renewcommand{\ge}{\geqslant}
\renewcommand{\Pr}{\mathbb{P}}

\newcommand{\Ex}{\mathbb{E}}

\title{Isoperimetry in product graphs}
\author{Sahar Diskin}
\author{Wojciech Samotij}
\address{School of Mathematical Sciences, Tel Aviv University, Tel Aviv 6997801, Israel}
\email{sahardiskin@mail.tau.ac.il, samotij@tauex.tau.ac.il}

\thanks{This research was supported by: the Israel Science Foundation grant 2110/22; the grant 2019679 from the United States--Israel Binational Science Foundation (BSF) and the United States National Science Foundation (NSF); and the ERC Consolidator Grant 101044123 (RandomHypGra).}

\begin{document}

\begin{abstract}
  In this short note, we establish an edge-isoperimetric inequality for arbitrary product graphs.
  Our inequality is sharp for subsets of many different sizes in every product graph.
  In particular, it implies that the $2^d$-element sets with smallest edge-boundary in the hypercube are subcubes and is only marginally weaker than the Bollobás--Leader edge-isoperimetric inequalities for grids and tori.
  Additionally, it improves two edge-isoperimetric inequalities for products of regular graphs proved by Erde, Kang, Krivelevich, and the first author and answers two questions about edge-isoperimetry in powers of regular graphs raised in their work.
\end{abstract}
\maketitle

\section{Introduction}
Given a graph $G$ with vertex set $V$, a key part of the (edge-)\emph{isoperimetric problem} is to determine, for every $k\in \mathbb{N}$, the quantity
\begin{align*}
  i_k(G)\coloneqq \min\left\{\frac{e_G(A,A^c)}{|A|} : A \subseteq V \wedge |A| = k\right\},
\end{align*}
where $e_G(A, A^c)$ is the number of edges of $G$ with exactly one endpoint in $A$.
For more details about discrete isoperimetric problems, we refer the interested reader to the surveys \cite{B94, B99, H04}.

In this note, we will consider the isoperimetric problem for \emph{product graphs}.
Instances of this problem have been studied in depth for several well-known product graphs, such as hypercubes~\cite{B67, H64, H76, L64}, Hamming graphs~\cite{L64}, grids, and tori \cite{BI91}.
Here, we will investigate the isoperimetric problem for arbitrary product graphs.
The motivation for considering this problem in such generality comes partially from the results of (and the questions posed in) the recent work~\cite{DEKK24}, where isoperimetric estimates played a crucial role in studying bond percolation on product graphs.

Given a positive integer $n$ and an arbitrary sequence of finite graphs $G_1, \dotsc, G_n$, the product graph $G_1\square\cdots\square G_n$ is the graph whose vertex set is $V(G_1)\times \dotsb \times V(G_n)$ and whose edges are all pairs $\{u,v\}$ for which there is an index $j\in \br{n}$ such that $u_jv_j \in E(G_j)$ and $u_m=v_m$ for all $m \neq j$.
In order to state our main result, we require the following definition.
Given an $m$-vertex graph $G$, let $\psi_G \colon [0, \log m] \to [0, \infty)$ be the \emph{convex minorant} of the function $\{\log k : k \in \br{m}\} \ni x\mapsto i_{e^x}(G)$; in other words, $\psi_G$ is the largest convex function satisfying $\psi_G(\log k) \le i_k(G)$ for all $k\in \br{m}$.\footnote{Here and throughout the paper, $\log$ denotes the natural logarithm.}
Observe that $\psi_G$ is piecewise linear and that the only points where its derivative is not continuous are of the form $\log k$ for some integer $k\in \br{m}$.
Further, $\psi_G$ is decreasing, as $i_k(G)\ge 0=i_{m}(G)$ for all $k\in \br{m}$.

\begin{thm}\label{th: main}
Let $n$ be a positive integer, let $G_1,\ldots, G_n$ be an arbitrary sequence of finite graphs, and let $\prodG\coloneqq G_1\square\cdots\square G_n$. For every $\emptyset \neq A\subseteq V(\prodG)$,
\[
    e_{\prodG}(A,A^c)\ge |A|\cdot \min\left\{\sum_{i=1}^n\psi_{G_i}(h_i) : 0\le h_i\le \log|V(G_i)|\wedge \sum_{i=1}^nh_i=\log|A|\right\}.
\]
In particular, if $G_1=\cdots=G_n=G$, then
\[
    e_{\prodG}(A,A^c)\ge |A|\cdot n\cdot \psi_G\bigl((\log |A|)/n\bigr).
\]
\end{thm}

Let us note that \Cref{th: main} gives a sharp bound for every $n$, every sequence $G_1,\dotsc, G_n$, and sets~$A$ of many different sizes.  To see this, assume first that $G_1=\cdots=G_n=G$ for some graph $G$ with $m$ vertices.
Consider arbitrary integers $k_1, k_2 \in \br{m}$ with $k_1 < k_2$ such that $\psi_G(\log k_i) = i_{k_i}(G)$ for both $i \in \br{2}$ and $\psi_G$ is linear on  $[\log k_1, \log k_2]$. Further, let $A_1, A_2 \subseteq V(G)$ be sets witnessing $|A_i| = k_i$ and $e_G(A_i, A_i^c) = i_{k_i}(G) \cdot k_i$ for both $i \in \br{2}$. Then, for all nonnegative integers $n_1$ and $n_2$ satisfying $n_1+n_2=n$, the set $A \coloneqq A_1^{n_1} \times A_2^{n_2} \subseteq V(G)^n$ satisfies
  \[
    \begin{split}
      \frac{e_{G^n}(A, A^c)}{|A|}
      & = n_1 \cdot i_{k_1}(G) + n_2 \cdot i_{k_2}(G)
        = n_1 \cdot \psi_G(\log k_1) + n_2 \cdot \psi_G(\log k_2) \\
      & = n \cdot \psi_G\left(\frac{n_1}{n} \cdot \log k_1 + \frac{n_2}{n} \cdot \log k_2\right)
        = n \cdot \psi_G\bigl((\log |A|) / n\bigr).
    \end{split}
\]
The above argument extends to product graphs that are not necessarily powers of a single graph.  In this general case, the lower bound on $e_{\prodG}(A, A^c)$ is achieved by sets $A$ of the form $A_1 \times \dotsb \times A_n$, where $A_i \subseteq V(G_i)$ satisfy $e_{G_i}(A_i, A_i^c) / |A_i|= i_{|A_i|}(G_i) = \psi_{G_i}(\log |A_i|)$ and, writing $\partial_- \psi(x)$ and $\partial_+ \psi(x)$ for the left and the right derivatives of $\psi$ at $x$,
\[
  \max_{i \in \br{n}} \partial_- \psi_{G_i}(\log |A_i|) \le \min_{i \in \br{n}} \partial_+ \psi_{G_i}(\log |A_i|),
\]
where we use the convention that $\partial_- \psi_{G_i}(0) = -\infty$ and $\partial_+ \psi_{G_i}(\log |V(G_i)|) = \infty$.  It is not hard to check that the above assumption on left and right derivatives ensures that the minimum in the statement of \Cref{th: main} is achieved at $(h_1, \dotsc, h_n) = (\log |A_1|, \dotsc, \log |A_n|)$.
\subsection*{Acknowledgement}

We thank Joshua Erde, Mihyun Kang, and Michael Krivelevich for their helpful comments and suggestions.

\subsection*{Organisation}

In \Cref{s: isoperimetry}, we present the (short) proof of \Cref{th: main}, and in \Cref{s: applications}, we discuss several applications of Theorem~\ref{th: main} and compare them with known results in the literature.

\section{Proof of Theorem \ref{th: main}}
\label{s: isoperimetry}

Our argument builds on the beautiful entropy-based proof of an optimal edge-isoperimetric inequality for the hypercube presented by Boucheron, Lugosi, and Massart in~\cite[Section 4.4]{BouLugMas13}.  
The \emph{entropy} of a discrete random variable $X$ taking values in a countable set $\mathcal{X}$ is the quantity $H(X)$ defined by
\[
  H(X) \coloneqq -\sum_{x\in\mathcal{X}} \Pr(X = x) \log \Pr(X = x).
\]
In particular, if $\mathcal{X}$ is finite and $X$ is uniform on $\mathcal{X}$, then $H(X) = \log|\mathcal{X}|$.  Further, given random variables $X$ and $Y$ taking values in countable sets $\mathcal{X}$ and $\mathcal{Y}$, respectively, we define the \emph{conditional entropy} of $X$ given $Y$, denoted $H(X \mid Y)$, to be the average entropy of the random variable $X$ conditioned on the outcome of $Y$; in other words,
\[
  H(X \mid Y) \coloneqq - \sum_{y \in \mathcal{Y}} \Pr(Y = y) \sum_{x \in \mathcal{X}} \Pr(X = x \mid Y = y) \log \Pr(X = x \mid Y = y).
\]

Let $n$ be a positive integer and suppose that $\prodG = G_1 \square \dotsb \square G_n$ for some arbitrary sequence $G_1, \dotsc, G_n$ of finite graphs.
Consider an arbitrary nonempty set $A \subseteq V(\prodG) = V(G_1) \times \dotsb \times V(G_n)$ and let $X = (X_1, \dotsc, X_n)$ be a uniformly chosen random vertex of $A$. For every $v\in V(\prodG)$ and each $i \in \br{n}$, denote by $v_{(i)}$ the projection of $v$ along the $i$th coordinate, that is, $v_{(i)}=(v_1,\ldots, v_{i-1}, v_{i+1}, \ldots, v_n)$. Further, given an $x \in A$, let $A_i(x) \subseteq V(G_i)$ denote the support of $X_i$ conditioned on $X_{(i)} = x_{(i)}$.
Our first key observation is that
\[
    e_{\prodG}(A,A^c)
    = \sum_{x \in A} \sum_{i=1}^n \frac{e_{G_i}\left(A_i(x), A_i(x)^c\right)}{|A_i(x)|}
    \ge \sum_{x \in A} \sum_{i=1}^n  i_{|A_i(x)|}(G_i).
\]
Denoting by $k_i$ the (random) size of $|A_i(X)|$, we may rewrite the above inequality as
\begin{equation}
    \label{eq:eAAc-probabilistic}
    e_{\prodG}(A, A^c) \ge |A| \cdot \sum_{i=1}^n \Ex[i_{k_i}(G_i)].
\end{equation}
By the definition of $\psi_{G_i}$ and by Jensen's inequality, we have, for every $i \in \br{n}$,
\[
    \Ex[i_{k_i}(G_i)] \ge \Ex[\psi_{G_i}(\log k_i)] \ge \psi_{G_i}\bigl(\Ex[\log k_i]\bigr).
\]
Our second key observation is that $\Ex[\log k_i]$ is precisely the conditional entropy $H(X_i \mid X_{(i)})$.
Substituting the above inequality into \eqref{eq:eAAc-probabilistic}, we conclude that
\[
    e_{\prodG}(A, A^c) \ge |A| \cdot \sum_{i=1}^n \psi_{G_i}\bigl(H(X_i \mid X_{(i)})\bigr).
\]
The main assertion of the theorem now follows as, for each $i \in \br{n}$, the function $\psi_{G_i}$ is decreasing, $0 \le H(X_i \mid X_{(i)}) \le H(X_i) \le \log |V(G_i)|$ for each $i$, and
\[
    \sum_{i=1}^n H(X_i \mid X_{(i)}) \le H(X) = \log |A|,
\]
by Han's inequality~\cite{Han78} (see \cite[Theorem~4.1]{BouLugMas13} for a compact statement).
Finally, if $G_1 = \dotsb = G_n = G$, then we may use the convexity of $\psi_G$ again to deduce that, for all sequences $(h_i)_{i=1}^n$ that sum to $\log |A|$,
\[
    \sum_{i=1}^n \psi_G(h_i) \ge n \cdot \psi_G\left(\sum_{i=1}^n \frac{h_i}{n}\right) = n \cdot \psi_G\left(\frac{\log|A|}{n}\right),
\]
as claimed. \qed

\vspace{-0.5em}
\section{Applications}\label{s: applications}

\subsection{Hamming graphs and the hypercube}

Let $K_m$ be the complete graph on $m$ vertices, so that $\prodG\coloneqq K_m^n$ is the Hamming graph $H(n,m)$.
Since $i_k(K_m) = m-k \ge (m-1) \cdot (1 - \log_m k)$ for all $k \in \br{m}$, where the inequality is equivalent to the inequality $(k-1)/(m-1) \le \log_m k$, which holds due to the concavity of $x \mapsto \log x$, we have
\begin{equation}
  \label{eq:psi of complete graphs}
  \psi_{K_m}(x) \ge (m-1) \cdot (1-x/\log m)
\end{equation}
for all $x\in [0,\log m]$.  Therefore, by \Cref{th: main}, for all nonempty $A \subseteq V(H(n,m))$, we have
\begin{equation}
  \label{eq:isoperimetry-complete-graphs}
  e_{H(n,m)}(A,A^c)\ge |A|\cdot (m-1)\left(n-\log_m|A|\right).
\end{equation}
Observe that~\eqref{eq:isoperimetry-complete-graphs} is sharp whenever $A$ induces a copy of $H(t,m)$ for some $t \in \br{n}$.
In this sense, one may view it as a weak version of the edge-isoperimetric inequality for Hamming graphs due to Lindsey~\cite{L64}.\footnote{Lindsey's inequality is the stronger statement that each initial interval in the lexicographic ordering of $\br{m}^n$ has the smallest edge-boundary among all sets of the same size.}
In particular, the case $m=2$, may be viewed as a weak version of the edge-isoperimetric inequality for the hypercube \cite{B67,H64,H76,L64}.

\subsection{The grid}
Let $P_m$ be the path with $m \ge 3$ vertices, so that $\prodG \coloneqq P_m^n$ is the $n$-dimensional $m \times \dotsb \times m$ grid.
Note that $i_k(P_m)=1/k$ for every $k \in \br{m-1}$ and that $i_m(P_m) = 0$.
For every $z \in [0, \log m)$, let $\ell_z$ be the line passing through the points $(z, e^{-z})$ and $(\log m, 0)$, that is, the line $y = e^{-z} \cdot (\log m - x)/(\log m - z)$.
Since the points $\{(\log k, 1/k) : k \in \br{m-1}\}$ lie on the graph of the convex function $x \mapsto e^{-x}$ and $\ell_{\log m -1}$ has the largest (that is, least negative) slope among all our lines $\ell_z$, we may deduce that
\begin{equation}
  \label{eq:Psi-Pm-lower}
  \psi_{P_m}(x) \ge
  \begin{cases}
    e^{-x} & \text{if $0 \le x \le \log m - 1$,}\\
    e/m \cdot (\log m - x) & \text{if $\log m - 1 \le x \le \log m$.}
  \end{cases}
\end{equation}
In fact, $\psi_{P_m}$ is the piecewise linear function defined by the points $(0,1), \dotsc, (\log k^*, 1/k^*)$, and $(\log m,0)$, where $k^* \in \br{m-1}$ is the index $k$ for which $\ell_{\log k}$ has the largest slope.  It is not hard to see that $k^*\in \left\{\lfloor m/e \rfloor, \lceil m/e\rceil\right\}$, but whether it is the floor or the ceiling of $m/e$ depends on the value of $m$.  For example, $k^*=\lfloor 3/e\rfloor =1$ when $m=3$, whereas $k^*=\lceil 5/e\rceil=2$ when $m=5$.

With the lower bound \eqref{eq:Psi-Pm-lower} in place, we can now use \Cref{th: main} to derive edge-isoperimetric inequalities for $\prodG$.  When $|A|\le (m/e)^n$, we have
\begin{align*}
    e_{\prodG}(A,A^c)\ge |A|\cdot n \cdot e^{-(\log|A|) / n} = n \cdot |A|^{1-1/n}
\end{align*}
and when $(m/e)^n \le |A|\le m^n/2$, we have
\begin{align*}
    e_{\prodG}(A,A^c)\ge |A|\cdot n \cdot \frac{e}{m}\bigl(\log m-(\log|A|)/n\bigr) = \frac{|A|}{m} \cdot e\log \frac{m^n}{|A|}.
\end{align*}
For comparison, Bollobás and Leader~\cite{BI91} showed that, for all $A \subseteq V(\prodG)$ with $|A|\le m^n/2$,
\[
  e_{\prodG}(A,A^c) \ge \frac{|A|}{m} \cdot \min\left\{r \cdot \left(\frac{m^n}{|A|}\right)^{1/r} : r\in \br{n}\right\}.
\]
Since the minimum above is achieved at $r=n$ whenever $|A| \le (m/e)^n$, our bound matches that of Bollobás and Leader in this range.  In the complementary range $(m/e)^n \le |A| \le m^n/2$, the ratio between the two bounds does not exceed
\[
  \sup\left\{\frac{\lceil \log x\rceil \cdot x^{1/\lceil \log x \rceil}}{e \log x} : 2 \le x \le e^n\right\} \le \sup\left\{\frac{e^{y-1}}{y} : 1/2 \le y \le 1\right\} = 2e^{-1/2} \le 1.214,
\]
where in the first inequality we substituted $y=\log x/\lceil \log x\rceil$ and used that $1/2 \le \log x / \lceil \log x \rceil \le 1$ for all $x \ge 2$.

\subsection{The torus}

Let $C_m$ be the cycle with $m$ vertices, so that $\prodG \coloneqq C_m^n$ is the $n$-dimensional discrete torus with side length $m$.
Since $i_k(C_m) = 2i_k(P_m)$ for all $k \in \br{m}$, we have $\psi_{C_m} = 2\psi_{P_m}$.  Thus, \Cref{th: main} and the estimate~\eqref{eq:Psi-Pm-lower} yield
\begin{align*}
  e_{\prodG}(A,A^c)\ge
  \begin{cases}
    2n \cdot |A|^{1-1/n} & \text{if $|A|\le (m/e)^n$},\\
    |A|/m \cdot 2e\log(m^n/|A|) &\text{if $|A| \ge (m/e)^n$}.
  \end{cases}
\end{align*}
For comparison, Bollobás and Leader~\cite{BI91} showed that, for all $A \subseteq V(\prodG)$ with $|A|\le m^n/2$,
\[
  e_{\prodG}(A,A^c) \ge \frac{|A|}{m} \cdot \min\left\{2r \left(\frac{m^n}{|A|}\right)^{1/r} : r\in \br{n}\right\},
\]
and hence, as in the case of grid graphs, our bound matches theirs whenever $|A| \le (m/e)^n$ and is off by a multiplicative factor of at most $2e^{-1/2}$ in the complementary range.

\subsection{Products of regular graphs}

For every $i\in \br{n}$, let $G_i$ be a $d_i$-regular graph on $m_i$ vertices, let $\prodG \coloneqq G_1\square\cdots\square G_n$, and note that $\prodG$ is also regular of degree $d \coloneqq d_1 + \dotsb + d_n$.  Note that the assumption that $G_i$ is $d_i$-regular implies that $i_k(G_i) \ge d_i - k + 1 = i_k(K_{d_i+1})$ for all $k \in \br{d_i+1}$; indeed, $\deg_{G_i}(v, A^c) \ge d_i - |A| + 1$ for all $v \in A \subseteq V(G_i)$.  Consequently, $\psi_{G_i}(x)\ge d_i \cdot (1-\log_{d_i+1}x)$ for all $x \in [0, \log m_i]$, see~\eqref{eq:psi of complete graphs}.  Thus, by \Cref{th: main}, 
\begin{align*}
    e_{\prodG}(A,A^c)&\ge |A|\cdot \left(d-\max\left\{\sum_{i=1}^n\frac{d_i\cdot h_i}{\log(d_i+1)} : 0\le h_i\le \log m_i\wedge \sum_{i=1}^nh_i=\log|A|\right\}\right)\\
    &\ge |A|\cdot \left(d-\max_{i \in \br{n}}\frac{d_i}{\log(d_i+1)} \cdot \log|A|\right) = |A|\cdot \bigl(d- D \cdot \log_{D+1}|A|\bigr),
\end{align*}
where $D \coloneqq \max_{i \in \br{n}} d_i$.  This substantially improves \cite[Theorem~1]{DEKK24}.

Assume further that each $G_i$ is connected, so that $i_k(G_i) \ge i_k(P_{m_i})$ for all $k \in \br{m_i}$.  It follows from~\eqref{eq:Psi-Pm-lower} that $\psi_{G_i}(x) \ge  e/m_i\cdot (\log m_i-x)$ for all $x \in [0, \log m_i]$.  Therefore, by \Cref{th: main},
\begin{align*}
    e_{\prodG}(A,A^c)&\ge |A|\cdot \min\left\{\sum_{i=1}^n\frac{e}{m_i} \cdot (\log m_i - h_i) : 0\le h_i\le \log m_i\wedge \sum_{i=1}^n h_i=\log|A|\right\}\\
    &\ge |A|\cdot \min\left\{ \sum_{i=1}^n \frac{e g_i}{m_i} : g_i \ge 0 \wedge \sum_{i=1}^n g_i = \log \frac{|V(\prodG)|}{|A|}\right\} = |A| \cdot \frac{e}{M} \cdot \log \frac{|V(\prodG)|}{|A|},
\end{align*}
where $M \coloneqq \max_{i \in \br{n}} m_i$.  When $M \ge 3$, this improves the respective lower bound on $e_{\prodG}(A, A^c)$  given by \cite[Theorem~2]{DEKK24} by a multiplicative factor of $e(1-1/M)\log M$.

\subsection{Powers of regular graphs}

Let $G$ be a connected $m$-vertex, $d$-regular graph and let $\prodG \coloneqq G^n$. For every $k\in \br{m-1}$, let $\ell_k$ be the line passing through $(\log k, i_k(G))$ and $(\log m, 0)$, that is, the line $y=i_k(G)\cdot(\log m-x)/(\log m-\log k)$.  Let $k^*$ be the smallest index $k$ such that $\ell_k$ has the least negative slope among all our lines and note that, for all $x\in [0,\log m]$,
\begin{equation}
  \label{eq:psi-G-regular}
  \psi_G(x)\ge i_{k^*}(G) \cdot \frac{\log m-x}{\log m-\log k^*}
\end{equation}
Let $y_G$ be the $y$-intercept of $\ell_{k^*}$. Note that $y_G \le d$ (as $i_1(G) = d$) and that $y_G=d$ if and only if $k^*=1$. Further, observe that $y_G = i_{k^*}(G) \cdot \log m / (\log m - \log k^*)$.
Hence, by \Cref{th: main}, 
\begin{equation}
  \label{eq:isoperimetry-product-regular}
  e_{\prodG}(A,A^c)\ge |A|\cdot \frac{i_{k^*}(G)}{\log m - \log k^*} \cdot \log \frac{m^n}{|A|} = |A|\cdot y_G\cdot (n-\log_m|A|).
\end{equation}
Since \eqref{eq:psi-G-regular} holds with equality for all $x \in [\log k^*, \log m]$, inequality~\eqref{eq:isoperimetry-product-regular} is tight for sets $A$ with many different sizes, see the construction described below the statement of \Cref{th: main}.

We now address two questions posed in~\cite{DEKK24}.
First, \cite[Question 7.1]{DEKK24} asked whether there are constants $c_G, C_G$ such that $i_{a}(\prodG)= c_G \cdot \log(m^n/a)+C_G$ for all $a\in \br{m^n}$.  In other words, \cite[Question 7.1]{DEKK24} asks whether $i_a(\prodG)$ is essentially linear in $\log a$.  The construction presented below the statement \Cref{th: main} shows that the lower bound on $i_a(\prodG)$ implied by the theorem is sharp whenever $\log a = (n_1/n) \cdot \log k_1 + (n_2/n) \cdot \log k_2$ for some $n_1, n_2$ satisfying $n_1 + n_2 = n$ and $k_1, k_2 \in \br{m}^2$ such that $[\log k_1, \log k_2]$ supports one of the linear pieces of $\psi_G$.  This fact implies that $i_a(\prodG)$ in not linear in $\log a$ whenever $\psi_G$ itself is not linear.  Since there are regular graphs $G$ for which $\psi_G$ has more than one linear piece (for example, when $G=C_m$ for $m\ge 5$), the answer to \cite[Question~7.1]{DEKK24} is negative.

Further, \cite[Question 7.2]{DEKK24} asked for a characterisation of $m$-vertex $d$-regular graphs $G$ for which sets of the form $B_t \coloneqq \{u\}^t\times V(G)^{n-t}$ have the smallest edge-boundary among all $m^{n-t}$-element sets of vertices of $\prodG$, for all $t \in \br{n}$.  We note that this is closely related to the classical problem of finding sufficient conditions for a graph to admit a nested sequence of sets that achieve the smallest edge-boundary (among all sets of a given size), see \cite{B99,H04} and references therein. Since
\[
  e_{\prodG}(B_t, B_t^c) = |B_t| \cdot t \cdot d = |B_t| \cdot \frac{d}{\log m} \cdot \log \frac{m^n}{|B_t|},
\]
it follows from~\eqref{eq:isoperimetry-product-regular} that a sufficient condition is $y_G=d$.  We will show below that, for large enough $n$, this is also a necessary condition.

Suppose that $G$ is an $m$-vertex, $d$-regular graph with $y_G < d$, let $k^* \in \{2, \dotsc, m-1\}$ be the index defined above, and let $S \subseteq V(G)$ be a $k^*$-element set witnessing $e_G(S, S^c) = |S| \cdot i_{k^*}(G)$.  Fix a small positive $\varepsilon$.  By Dirichlet's approximation theorem, there exist positive integers $s$ and $t$ such that
\begin{equation}
  \label{eq:Dirichlet}
  \bigl| s \log m - t \log (m/k^*) \bigr| \le \varepsilon/2,
\end{equation}
which implies that $(1-\varepsilon)m^t \le (k^*)^t \cdot m^s \le (1+\varepsilon)m^t$.  Consider the graph $\prodG \coloneqq G^{s+t}$ and sets of vertices $A \coloneqq S^t \times V(G)^s$ and $B \coloneqq \{u\}^s \times V(G)^t$.  Note that, by~\eqref{eq:Dirichlet},
\[
  e_{\prodG}(A, A^c) = |A| \cdot t \cdot i_{k^*}(G) = |A| \cdot \frac{y_G}{\log m} \cdot t \log(m/k^*) \le |A| \cdot y_G \cdot (s + \varepsilon) \le (1+2\varepsilon) \cdot m^t \cdot y_G\cdot s.
\]
Let $C$ be a set of size exactly $m^t$ that is obtained by adding to / removing from $A$ at most $\varepsilon m^t$ vertices in an arbitrary manner.  Since $\Delta(\prodG) = (s+t)d$, we clearly have
\[
  e_{\prodG}(C, C^c) - e_{\prodG}(A, A^c) \le \varepsilon m^t \cdot (s+t)d \le \varepsilon m^t s d \left(1 + \frac{\log m + \varepsilon/2}{\log(m/k^*)}\right),
\]
where the second inequality follows from~\eqref{eq:Dirichlet}.
Since we assumed that $y_G < d$, it is clear that choosing $\varepsilon$ sufficiently small (as a function of $m$ and $d-y_G$ only) gives $e_{\prodG}(C, C^c) < m^t s d = e_{\prodG}(B, B^c)$. This means that the set $B$ does not have the smallest edge boundary among all sets of $m^t$ vertices of $\prodG$.

\bibliography{isoperimetry}
\bibliographystyle{abbrv}

\end{document}